\documentclass[11pt,a4paper]{article}
\usepackage[OT1]{fontenc}
\usepackage[latin1]{inputenc}
\usepackage[english]{babel}
\usepackage{amsfonts}
\usepackage{amsthm}
\def\grass{{\cal P}(M_1)}

\newtheorem{teo}{Theorem}[section]
\newtheorem{prop}[teo]{Proposition}
\newtheorem{lem}[teo]{Lemma}
\newtheorem{coro}[teo]{Corollary}
\newtheorem{defi}[teo]{Definition}
\theoremstyle{definition}
\newtheorem{rem}[teo]{Remark}

\newcommand{\p}{{\sf \,p\,}}
\newcommand{\h}{{\cal H}}
\newcommand{\um}{ {\cal O}(\p)  }

\newcommand{\pei}{{\left<\right.}}
\newcommand{\ped}{{\left>\right.}}

\begin{document}

\title{\vspace*{0cm}Weak Riemannian manifolds from finite index
subfactors\footnote{2000 MSC. Primary 58B20;  Secondary 46L10, 53C30, 53C22.}}
\date{}
\author{Esteban Andruchow and Gabriel Larotonda}

\maketitle

\abstract{\footnotesize{\noindent Let $N\subset M$ be a finite  Jones' index
inclusion of II$_1$ factors, and denote by $U_N\subset U_M$ their unitary
groups. In this paper we study the homogeneous space $U_M/U_N$, which is a
(infinite dimensional) differentiable manifold, diffeomorphic to the orbit
$$
{\cal O}(p) =\{u p u^*: u\in U_M\}
$$
of the Jones projection $p$ of the inclusion. We endow ${\cal O}(p) $ with a Riemannian
metric, by means of the trace on each tangent space. These are  pre-Hilbert
spaces (the tangent spaces are not complete), therefore ${\cal O}(p)$ is a weak Riemannian manifold. We show that ${\cal O}(p)$ enjoys certain properties similar to classic Hilbert-Riemann manifolds. Among them, metric
completeness of the geodesic distance, uniqueness of geodesics of the
Levi-Civita connection as minimal curves, and  partial  results on the existence
of minimal geodesics. For instance, around each point $p_1$ of ${\cal O}(p)$, there is a
ball $\{q\in {\cal O}(p):\|q-p_1\|<r\}$ (of uniform radius $r$) of the usual norm of $M$, such that any point $p_2$ in the ball is joined to
$p_1$ by a unique geodesic, which is shorter than any other piecewise smooth
curve lying inside this ball. We also give an intrinsic (algebraic) characterization of the directions of degeneracy of the submanifold inclusion ${\cal O}(p)\subset {\cal P}(M_1)$, where the last set denotes the Grassmann manifold of the von Neumann algebra generated by $M$ and $p$.}\footnote{{\bf Keywords and
phrases:} homogeneous space, short geodesic, Levi-Civita connection, Riemannian submanifold, totally geodesic submanifold, finite index inclusion, von Neumann $II_1$ subfactor, Jones' projection, trace quadratic norm.}}

\setlength{\parindent}{0cm} 

\section{Introduction}
Let $N\subset M$ be a finite index inclusion of II$_1$ factors,
and let $U_N\subset U_M$ be their  unitary groups. In this paper
we study the homogeneous space $U_M/U_N$ as an example of a weak Riemannian manifold, i.e. an infinite dimensional manifold with a Riemannian metric which makes the tangent spaces pre-Hilbert spaces. This paper is a continuation of \cite{as}, where the
basic topological and differential facts of this space were
established. Mainly, that there exists a concrete model for this
 homogeneous space, which is diffeomorphic to it. Let us
describe this model. Consider Jones' basic construction: let $\p$
denote the Jones projection of the inclusion, i.e. the orthogonal
projection
$$
\p:L^2(M,\tau)\to L^2(N,\tau),
$$
where $L^2(M,\tau)$ and  $L^2(N,\tau)$ are the GNS Hilbert spaces
of the trace $\tau$ (we denote indistinctly by $\tau$ the traces
of $N$ and $M$). Consider the von Neumann algebra $M_1$ of
operators in $L^2(M,\tau)$ generated by $M$ and $\p$. Then
$$
U_M/U_N\simeq \um:=\{u\p u^*: u\in U_M\}\subset M_1.
$$
This natural bijection is a homeomorphism. In \cite{as}
it was shown that this orbit $\um$ is a submanifold of the
Grassmann manifold $\grass$ (=selfadjoint projections) of $M_1$. Moreover, this
regularity  property in fact is equivalent to the finite index
condition. The Grassmann manifold  of a $C^*$-algebra is a well
behaved Finsler manifold, with a linear connection coming from a reductive
structure, and the Finsler metric induced by the
usual norm of the algebra \cite{cpr1,pr}. The geodesics of the linear
connection are minimal curves between points (=projections) lying at norm
distance
less than $1$.

It is well known \cite{jones,popa} that $M_1$ is also a
type II$_1$ factor, and that the inclusion $M\subset M_1$ has the
same index as the original inclusion. This gives rise the the so called {\it
basic construction}, which we recall later.  Denote also by $\tau$ be the
trace of $M_1$ (which by restriction gives the traces of $M$ and
$N$).

In von Neumann algebras the relevant topologies are the weak topologies. The
$\sigma$-strong, strong and weak operator topologies coincide in $\grass$. This
set is a manifold only in the norm topology. Therefore, whenenever we use the
word smooth, applied to a curve or a map, it is meant that is smooth in the norm
structure.  The weak topologies in  $\grass$ are metrized by the 2-norm induced
by the trace,
$\|x\|_2=\tau(x^*x)^{1/2}$. Although $\grass$  is {\bf not} a submanifold
of the GNS Hilbert space completion $L^2(M_1,\tau)$, it is a complete metric
space, the geodesic metric is equivalent to the $2$-metric, and the geodesics
of the linear connection mentioned above behave nicely with
respect to the 2-norm (in fact, with the $k$-norms, for $k\ge 2$): if
two projections lie at (norm) distance less than $1$, then the geodesic
that joins them is shorter than any other piecewise smooth curve,
when the length is measured with the 2-norm
\cite{argrassmannians}.

In this paper we study the metric induced in $\um$ by the trace
inner product $<x,y>=\tau(y^*x)$ at each tangent space. The
tangent spaces
$$
(T\um)_q=\{xq-qx: x\in M, x^*=-x\}
$$
are not complete with this metric. Therefore this study does not
fit in the classical (infinite dimensional) Riemannian theory
\cite{lang}. Certain computations, though, can be carried out nicely, as in
the classical case: e.g. the geodesics of the Levi-Civita connection
can be explicitly computed. However, the fact that the tangent
spaces are not complete allow for certain unusual phenomena. For
instance, one cannot find normal neighborhoods around each point.
Nevertheless certain facts do hold. Let us remark, for example, that $\um$,
with the geodesic distance induced by this incomplete metric, is
in fact a complete metric space (Theorem \ref{completo}). Also we prove
(Theorem \ref{unici}), that if there exists a curve that is piecewise
smooth  (in the norm structure), such that it has minimal
length for the trace  metric, then it is a geodesic
of the linear connection. As for existence of minimal geodesics, we show
(Theorem 5.3) that there exists a radius $r$, such that if
$\|p_1-p_2\|<r$, $p_1,p_2\in \um$, then there exists a geodesic of the
Levi-Civita connection of the trace, with $p_1$ and $p_2$ as its endpoints,
which is shorter than any other piecewise smooth curve lying inside the ball of
center $p_1$ and radius $r$. Regarding the submanifold inclusion $\um\subset {\cal P}(M_1)$ (see \cite{as}), we show (Theorem 6.3) that the Riemannian curvature of this inclusion is related to the inclusion $\{ x^*=-x \in \ker E:\; x^2\in N\}\subset M$, and this set is exactly the set of directions of degeneracy in $\um$ (\textit{i.e.} the directions $x$ in the tangent bundle such that the geodesic of ${\cal P}(M_1)$ with initial speed $x$ is also a geodesic of $\um$).

\medskip

The contents of the paper are as follows. In Section 2 we recall certain basic
facts, as the Jones basic construction, and establish topologic properties of
$\um$, both in the norm and weak topologies. In Section 3 we introduce the
Riemannian structure in $\um$ by means of the trace $\tau$. We show that though
the tangent spaces are incomplete, the metric space $\um$ with the Riemannian
(or geodesic) distance is complete. Also we study the horizontal lifting of
curves, and its metric properties. Based on a result on the local convex
structure of the unitary group in the $k$-norms proved in \cite{ar}, we prove
that geodesics are short curves (in the $2$-metric) among curves which do not
exceed certain length (measured in the usual norm Finsler metric). This result
is the key fact to our minimality results in Section 5, where we also prove the
uniqueness of geodesics as minimal curves. In Section 6 we study the geodesics of the Riemannian inclusion $\um\subset {\cal P}(M_1)$.

\section{Topological considerations}

Let $N\subset M$ be an inclusion of II$_1$ factors with trace $\tau$, and finite
index $[M:N]=\lambda^{-1}$. Let $L^2(M,\tau)$ and $L^2(N,\tau)$ be the GNS
Hilbert spaces of $\tau$ and let
$\p: L^2(M,\tau)\to L^2(N,\tau)$ be the (Jones) orthogonal projector. This
projection induces the
unique trace-preserving conditional expectation $E:M\to N$. Let $M_1$ be the von
Neumann algebra of operators in $L^2(M,\tau)$ generated by $\p$ and $M$.
This construction, known as the basic construction, enjoys several properties.
Among them \cite{jones,popa}:
\begin{enumerate}
\item
$M_1$ is a finite factor, $[M_1:M]=[M:N]=\lambda^{-1}$.
\item
 $\p x \p=E(x)\p$ for any $x\in M$.
\item
 $\{\p\}'\cap M=N$.
 \item
 $N\ni x\mapsto x\p \in N\p= \p M_1 \p$ is a $*$-isomorphism.
  \item
 $M_1 \p=M\p$, and therefore $M\p$ is a closed linear subspace of $M_1$.
 \item
 The  map $M\ni a\mapsto a\p \in M\p$ is a linear isomorphism, with $\|a\|\ge
\|a\p\|\ge \sqrt{\lambda}\|a\|$.
 \item
 If we denote by $E_1:M_1\to M$ the unique $\tau$-preserving conditional
expectation, then $E_1(\p )=\lambda$.

 \item
 $E(x^*x)\ge \lambda x^*x$ for any $x\in M$.
 \end{enumerate}

\bigskip

\begin{defi}\label{defin}
Let $\um=\{u\p u^*:u\in U_M\}$ be the $U_M$-unitary orbit of $\p$, viewed as a
subset
of $M_1$. Since $\{\p\}'\cap M=N$, $\um$ can be identified with the quotient
$U_M/U_N$ via the
map $\ell_{\p}:U_M\to M_1$,  $\ell_{\p}(u)= u\p u^*$.
\end{defi}

Note that $\um\simeq U_M/U_N$ is a topological isomorphism (for the norm induced
topologies), and when the index is finite, it can be proved that $\ell_{\p}$ is
a fibration (see \cite{as}), and admits smooth cross-sections. In that paper it
was proved that $\um$ is a smooth submanifold of ${\cal P}(M_1)$ (the set of
projections of $M_1$) if and only if the index of the inclusion $N\subset M$ is
finite. The following result is certainley well known, and shall be used later. Using this fact below, we may construct local cross sections for $\um$, from local cross sections for $\grass$.

\begin{lem}\label{lemaseccion}
Suppose that $\omega \in U_{M_1}$ satisfies that $\omega \p \omega^*\in\um$.
Then there exists $u\in U_M$ such that $u\p=\omega\p$. Moreover, if $\omega(t)$
is a smooth curve of unitaries in $M_1$ with the same property, then the curve
$u(t)$ of unitaries in $M$ can also be chosen smooth.
\end{lem}
\begin{proof}
By the facts 5 and 6 in the list of properties of the basic construction listed
at the beginning of Section 2, there exists $m\in M$ such that $m\p=\omega \p$.
We claim that $m$ is a unitary element.
Note that $\omega \p\omega^*=\omega \p(\omega\p)^*= m\p m^*\in \um$. Therefore
there exists a unitary element $v\in U_M$ such that $m\p m^*=v\p v^*$. Then
$v^*m \p m^*v=\p$. If we multiply by $\p$ on both sides, we obtain
$$
\p=\p v^*m \p m^*v\p=E(v^*m)\p E(m^*v)\p=E(v^*m)E(m^*v)\p,
$$
because $\p$ commutes with $N$, and then by further properties of the basic
construction,
$$
E(v^*m) E(m^*v)=E(v^*m)\p E(v^*m)^*\p=1.
$$
Since $N$ is finite, this implies that $E(v^*m)$ is a unitary element. Then
$\p=v^*m\p m^*v\p=v^*m E(m^*v)\p$, or equivalently, using that $E(v^*m)$ is
unitary,
$$
v^*m \p=E(v^*m)\p.
$$
By property 6 of the list, this implies that $v^*m =E(v^*m)$, and therefore
$m\in U_M$.
If $\omega(t)$ is a smooth curve in $U_{M_1}$ with
$\omega(t)\p\omega^*(t)\in\um$, one proves that $u(t)\in U_M$ is also smooth,
because $m$ above is obtained by composition of the map $M_1\ni\omega\mapsto
\omega \p\in M\p$ with the inverse linear isomorphism of property 6, $M\ni
a\mapsto a\p\in M\p$.
\end{proof} 

This lemma above shows that  $u$ can be
computed by
$$
u=\frac{1}{\lambda}E_1(w\p).
$$

It is known \cite{argrassmannians} that the Grassmann manifold $\grass$ has geodesical radius
 $1$: two projections $p_1$ and $p_2$ in $M_1$ such that
$\|p_1-p_2\|<1$ are joined by a unique curve $\gamma(t)=e^{tx}p_1e^{-tx}$,
with $x^*=-x$ in $M_1$, $x$ $p_1$-codiagonal, and $\|x\|< \pi/2$. The
condition that $x$ is $p_1$-codiagonal describes precisely that $\gamma$ is a
geodesic of the connection in $\grass$, and means that
$p_1xp_1=(1-p_1)x(1-p_1)=0$, or equivalently, that $x=xp_1+p_1x$. In particular
one may define a distinguished  {\it exponential } local cross section for the map
$$
U_{M_1}\to \grass , \ w\mapsto wp_1w^*,
$$
namely
$$
s_{p_1}:\{p_2\in\grass: \|p_2-p_1\|<1\} \to U_{M_1} , \ \ s_{p_1}(p_2)=e^x.
$$
Putting these two facts together enables us to obtain that if $q\in\um$ such that $\|q-\p \|<1$,
then
$$
\theta_\p(q)=\frac{1}{\lambda}E_1(s_\p(q)\p),
$$
defines a real analytic (in the norm structure) local cross section for the homogeneous space $\um$:
$$
\theta_\p(q)\in U_M \ \ \ \hbox{ and } \ \ \ \theta_\p(q)\p \theta_\p(q)^*=q.
$$

Our next result states that the identification $\um\simeq U_M/U_N$ also works in
the strong operator topology. Note that since $N\subset M$ are finite, both
$U_M$, $U_N$ are complete topological groups in the strong operator topology,
which is metrized by the $2$-norm $\|\ \|_2$.
To prove our statement, we shall  need the following result, which is certainly
not unknown to specialists (see  Lemma 5.3 in \cite{asnilpotents} for a proof):
\begin{lem}\label{lemal2}
Let $N$ be a finite von Neumann algebra and let $a_n\in N$ such that $\|a_n\|\le
1$ and $a_n^*a_n\to 1$ strongly. Then there exist unitaries $v_n\in N$ such that
$a_n-v_n\to 0$ strongly.
\end{lem}
\begin{prop}\label{biyeccion}
The natural bijection
$$
U_M/U_N \to \um , \ \ [u]\mapsto u\p u^*
$$
is a homeomorphism with the topologies  induced by the strong
operator topology.
\end{prop}
\begin{proof}
Suppose that $\{[u_d]\}_{d\in D}$ is a net in the quotient $U_M/U_N$ which
converges to $[u]$ in the topology induced by the strong operator topology. This
implies that there exist unitaries $v_\alpha\in U_N$ such that $u_\alpha
v_\alpha\to u$ strongly. Since $M$ is finite, this implies that also
$v_\alpha^*u_\alpha^* \to u^*$ strongly. Then, using that multiplication is also
strongly continuous,
$$
u_\alpha \p u_\alpha^*=u_\alpha v_\alpha\p v_\alpha^*u_\alpha^*\to u\p u^*.
$$
Note that $\um$ is a bounded set in $M_1$, and therefore the strong topology is
metrizable with the $2$-norm. In particular we may replace all nets by sequences
in our arguments. By the strong continuity of the action of $U_M$ on $\um$
(again using that $M$ is finite), we only need to prove continuity of the
inverse map at $\p$. Suppose that $u_n\p u_n^* \to \p$. Then we claim that
$\|u_n E(u_n^*)-1\|_2\to 0$. Indeed
\begin{eqnarray}
\|u_n E(u_n^*)-1\|_2^2 & = & \tau(E(u_n) E(u_n^*))+1-\tau(u_n
E(u_n^*))-\tau(E(u_n) u_n^*)\nonumber \\
&= & 1-\tau(E(u_n) E(u_n^*)),\nonumber
\end{eqnarray}
because $\tau(u_n E(u_n^*))=\tau(E(u_n E(u_n^*)))=\tau(E(u_n)
E(u_n^*))$ and similarly for the other term. Therefore it suffices to show
that $\tau(E(u_n) E(u_n^*))\to 1$ strongly.
Note that $u_n\p u_n^*\to \p$ strongly implies that
$$
\p u_n\p u_n^*\p=E(u_n)E(u_n^*)\p\to \p
$$
strongly. Note that $E(u_n)E(u_n^*)\in N$ and recall that the map $N\to N\p$
is a $*$-isomorphism, which is clearly normal. Therefore $E(u_n)E(u_n^*)\to
1$ strongly, and thus $\tau(E(u_n)E(u_n^*))\to 1$.
We may apply Lemma \ref{lemal2} to the elements $a_n=E(u_n^*)$. It follows
that there exist unitaries $v_n$ in $N$ such that $E(u^*_n)-v_n\to 0$
strongly, and therefore
$u_nv_n-u_nE(u_n^*)\to 0$ strongly. Thus
$$
u_nv_n\to 1,
$$ which completes the proof.
\end{proof}

\begin{prop}\label{completitud}
The orbit $\um$ is a complete metric space in the $2$-norm.
\end{prop}
\begin{proof}
Suppose that $\{p_n\}$ is a Cauchy sequence in $\um$, with $p_n=u_n\p u_n^*$,
$u_n\in U_M$. Then $p_n$ converges strongly to an operator $q$ acting in
$L^2(M_1,\tau)$. Indeed, for each $\eta\in L^2(M_1,\tau)$, $\{p_n\eta\}$ is a
Cauchy sequence in $L^2(M_1,\tau)$. Let us prove this fact.
If $\eta=x\in M\subset L^2(M_1,\tau)$,
\begin{eqnarray}
\|p_nx-p_kx\|_2^2 & = & \tau(x^*(p_n-p_k)^2 x)=\tau((p_n-p_k)x^*x(p_n-p_k)) \nonumber\\
 & \le & \|x\|^2\tau((p_n-p_k)^2) =\|x\|^2\|p_n-p_k\|_2^2.\nonumber
\end{eqnarray}
In general, there exists $x\in M$ such that $\|x-\eta\|_2<\epsilon/2$. And
therefore
\begin{eqnarray}
\|(p_n-p_k)\eta\|_2 &\le & \|(p_n-p_k)x\|_2+\|(p_n-p_k)(\eta-x)\|_2 \nonumber\\
 &\le & \|(p_n-p_k)x\|_2+2\|\eta-x\|_2 <\|(p_n-p_k)x\|_2+\epsilon.\nonumber
\end{eqnarray}
Therefore $p_n$ converges strongly to a linear operator in $\h$, which is
bounded by the uniform boundedness principle, and lies in $M_1$, i.e. $p_n\to
q\in M_1$. By strong
continuity of the product and the adjoint ($M_1$ is  finite), clearly
$q^2=q^*=q$. Moreover, since $\tau$ is normal and $\tau(p_n)=\tau(\p)=\lambda$,
it follows that $\tau(q)=\lambda$. Therefore $q=w\p w^*$ for $w\in U_{M_1}$, and
thus by the properties of the basic construction, there exists $m\in M$ with
$E(m^*m)=1$ such that $q=m\p m^*$. Using the conditional expectation
$E_{M_1}:M_1\to M$, which verifies that $E_{M_1}(\p)=\lambda$, one obtains that
$$
E_{M_1}(p_n)=u_nE_{M_1}(\p)u_n^*=\lambda,
$$
and therefore $E(q)=\lambda$. Then
$$
\lambda=E_{M_1}(m\p m^*)=m E_{M_1}(\p)m^*=\lambda m m^*,
$$
that is, $m\in U_M$ and therefore $q\in\um$.
\end{proof}

In a remarkable paper \cite{pota}, S. Popa and M. Takesaki proved
the following result \cite[Lemma 3]{pota}, based on the theory of
continuous selections by E. Michael \cite{michael}. Suppose that the  separable
II$_1$ factor $N$ has the following property, which for brevity we will refer as
the {\it scaling } property:  the tensor product $N\otimes {\cal B}(H)$ ($H$ a
separable Hilbert space) admits a one parameter automorphism group $\{\theta_s:
s\in \mathbb R\}$ scaling the trace of $N\otimes {\cal B}(H)$, i.e.
$\tau\circ\theta_s=e^{-s} \tau$, $s\in \mathbb R$, with $\tau$ a faithful
semi-finite normal trace in $N\otimes {\cal B}(H)$. Then $U_N$ is contractible
in the strong operator topology, and moreover, if $M$ is another factor with
$N\subset M$, then the quotient map
$$
U_M\to U_M/U_N
$$
admits a global cross section when these spaces are considered with the
(topologies induced by) the strong operator topology. The family of algebras
enjoying these scaling properties includes the hyperfinite factor ${\cal R}$.
A straightforward consequence of this result follows:
\begin{prop}
Suppose that $N$ has the scaling property. Then the map
$$
U_M \to \um , \ \ u\mapsto u\p u^*
$$
is a fibre bundle, with the topologies induced by the strong
operator topology.
If moreover also $M$ has the scaling property, then $\um$ is contractible in the
strong operator topology.
\end{prop}
\begin{proof}
The proof follows using the homotopy exact sequence of the bundle above, noting
that $\um\simeq U_M/U_N$ in the strong topology as well.
\end{proof}

\bigskip

Let us finish this section with the following result, concerning the homotopic
structure of $\um$ in the {\it norm} topology, for general $N\subset M$.

\begin{prop}
The space $U_M/U_N$ (hence also $\um$) is simply connected in the norm induced
topology.
\end{prop}
\begin{proof}
In the norm topology, the map
$$
\ell_\p: U_M\to \um , \ \ \ell_\p(u)=u\p u^*
$$
is a fibre bundle with fibre $U_N$. It therefore induces an exact
 sequence for the homotopy groups
$$
\dots\to \pi_1(U_N)\mathop\to^{T} \pi_1(U_M)\to \pi_1(\um)\to
\pi_0(U_N)\to\cdots
$$
H. Araki, M.B. Smith and L. Smith proved in \cite{arass}
that $\pi_1(U_N)=\pi_1(U_M)\simeq \mathbb R$. Moreover, they showed that this
isomorphism is given as follows. Any given loop $\beta$ in $U_M$, with base
point $1$, is homotopic to a concatenation of loops of the form
$\alpha_p(t)=e^{2\pi itp}$, for $p$ a projection. In other words, these loops
$\alpha_p$ generate $\pi_1(U_M)$. The isomorphism is given by the map which
sends the homotopy class of $\alpha_p$ to the real number $\tau(p)$. Therefore,
via this identification, the map $T$ above, restricted to this set of generators
is given by the inclusion map
$$
\{\tau(p): p\in {\cal P}(N)\}\hookrightarrow \{\tau(q): q\in {\cal P}(M)\},
$$
which is surjective, because both set equal the unit interval.
Hence $\pi_1(\um)=0$.
\end{proof}

\section{Riemannian structure of $\um$}

Using the $2$-norm $\|\ \|_2$ we can measure the length of curves in the
standard fashion:
$$
L_2(\alpha)=\int_0^1\|\dot{\alpha}(t)\|_2 dt.
$$
Let $M_{ah}=\{x\in M:x^*=-x\}$, $N_{ah}=\{x\in N:x^*=-x\}$. These spaces
identify with the Banach-Lie algebras of $U_M$ and $U_N$ respectively (i.e.
$U_M=\exp(M_{ah})$ and $U_N=\exp(N_{ah})$). Using that the map
$$
\ell_q:U_M\to \um , \ \ \ell_q(u)=u q u^*
$$
is a C$^\infty$-submersion for $q\in\um$, one can compute the tangent spaces of
$\um$ at any $q$. The differential of $\ell_q$ at $1$ is given by
$$
\delta_q(x)=d(\ell_q)_1(x)=xq-qx, \ \ x\in M_{ah}.
$$
Therefore
$$
(T\um)_q=\{xq-qx:x\in M_{ah}\}.
$$

We have an identification \textit{via} the differential of the quotient map
$U_M\mapsto U_M/U_N$
$$
(T U_M/U_N)_{\p}\simeq M_{ah}/N_{ah}
$$

\begin{rem}
Since $\tau\circ E=\tau$, if we put $\h=\ker E\cap M_{ah}$ we have
$$
M_{ah}=N_{ah}\oplus \h
$$
where $\h$ acts as an orthogonal supplement of $N_{ah}$ in $M_{ah}$ (relative to
the inner product given by $\tau$). Note that $\h$ is $Ad(U_N)$-invariant,
namely $u\h u^*=\h$ for any $u\in U_N$; let $\h_q=u\h u^*=Ad_{u}(\h)$ $(q=u\p
u^*, \,u\in U_M)$.
Likewise, we will write $E_q$ to denote the translated conditional expectation
$E_q=Ad_u\circ E\circ Ad_{u^*}$, which maps $M$ onto the von Neumann algebra
$Ad_u(N)$.

Therefore the maps $\delta_q$ give isomorphisms:
$$
\delta_q:\h_q\to(T\um)_q.
$$
\end{rem}

\begin{defi}
Let $\kappa_q$ be the inverse of this linear map, namely $\kappa_q(z q-q
z)=z$ for $z\in \h_q$, $\kappa_q: (T \um)_q\to \h_q$. Equivalently,
$\kappa_q(y)=z$ where $z$ is the only element of $\h_q=\ker E_q$ such that
$\delta_q(z)=y$. We will omit the subindex when $q=p$, i.e.
$\kappa_{\p}=\kappa$,
$\delta_{\p}=\delta$.
\end{defi}
Note that $\kappa_q$ is the $Ad-$translation of $\kappa$ to the point $q=u\p
u^*\in \um$, namely
$$
\kappa_q(v)=Ad_u\circ \kappa\circ Ad_{u^*}(v)=u \kappa(u^*vu)u^*.
$$

\begin{lem}\label{qiso}
The isomorphism $\delta_q: T_{\p}\um\simeq\h$ above is almost isometric,
\it{i.e}  if $q\in\um$,
$$
\|\delta_q(z)\|_2=\sqrt{2\lambda}\|z\|_2
$$
for any $z\in \h_q$.
\end{lem}
\begin{proof}
We use the properties of $E$ and $\tau$ in a II$_1$ factor with finite index
$\lambda$:
$$
\|z\p-\p z\|_2^2=2\tau(\p z\p z-\p
z^2\p)=-2\tau(E(z^2)\p)=-2\lambda\tau(E(z^2))=-2\lambda\tau(z^2)=2\lambda\|z\|_2^2
$$
The identity now follows from the unitary invariance of the 2-norm.
\end{proof}

In order to compute the linear connection induced by the trace metric, one needs
compute first the orthogonal projection from the ambient space, i.e. $M_{1\
h}=\{x\in M_1:x^*=x\}$ onto the tangent spaces of $\um$. Recall from the
properties of the basic construction, that $M_1\p=M\p$ induces a map
$$
R: M_1\to M, \ \ R(x)=m
$$
characterized as the unique element $m\in M$ such that $m\p=x\p$. This map can
be computed in terms of $E_1$: $E_1(x\p)=E_1(m\p)=m\lambda$, i.e.
$$
R(x)=\frac{1}{\lambda}E_1(x\p).
$$
Then we claim that the $\tau$-orthogonal projection
$$
\Pi_\p:M_{1\ h}\to (T\um)_\p=\{z\p-\p z:z\in M_{ah}\}
$$
is given by
$$
\Pi_\p(x)=\frac12[R(x)-R(x)^*,\p],
$$
or equivalently
\begin{equation}\label{proyeccionaltangente}
\Pi_\p(x)=\frac{1}{2\lambda}[E_1(x\p-\p x),\p],
\end{equation}
where $[\ , \ ]$ is the usual commutator of operators.
First note that $\Pi_\p$ projects onto $(T\um)_\p$. Apparently it takes values in this
linear subspace. If $z\p-\p z\in (T\um)_\p$ ($z\in \h$), then
$(z\p-\p z)\p=z\p$ because $\p z\p=0$. Therefore $R(z\p-\p z)=z$ and
$$
\Pi_\p(z\p-\p z)=\frac12[z-z^*,\p]=z\p-\p z.
$$
Finally, note that it is symmetric.
Indeed, if $x,y\in M_{1 \ h}$,
\begin{eqnarray}
<\Pi_\p(x),y> & = & \frac{1}{2\lambda}\tau([E_1(x\p-\p x),\p]y) \nonumber \\
&= & \frac{1}{2\lambda}\{\tau(E_1(x\p)\p y)-\tau(E_1(\p x)\p
y)+\tau(E_1(x\p)y\p)+\tau(E_1(\p x)y\p)\}\nonumber\\
&= &\frac{1}{2\lambda}\{\tau(E_1(x\p)E_1(\p y))-\tau(E_1(\p x)E_1(\p
y)) \nonumber\\
&&+  \tau(E_1(x\p)E_1(y\p))+\tau(E_1(\p x)E_1(y\p))\},\nonumber
\end{eqnarray}
which is clearly a symmetric expression in $x$ and $y$.
In order to obtain the symmetric projections onto the other tangent spaces
$(T\um)_q$, $q\in\um$ one translates $\Pi_\p$ covariantly via the action of
$U_M$, namely if $q=u\p u^*$ for $u\in U_M$
$$
\Pi_q:M_{1 \ h}\to (T\um)_q, \ \ \Pi_q=Ad(u)\circ\Pi_\p\circ Ad(u^*),
$$
which does not depend on the choice of $u$ by general reasons. Nevertheless note
that
$$
\Pi_q(x)=\frac{1}{2\lambda}u[E_1(u^*xu\p-\p
u^*xu),\p]u^*=\frac{1}{2\lambda}[E_1(xq-qx),q].
$$

\begin{rem}
Let $X$ be a smooth vector field along a curve $\gamma\subset\um$, then the
covariant derivative induced by the trace inner product is given by
$$
\frac{D X}{dt}=\Pi_\gamma(\dot{X}).
$$
We shall call this connection the Levi-Civita connection of $\um$ for obvious
reasons.

This  derivative  is compatible with the metric and has no torsion, these facts
follow from formal considerations (as in classical Riemannian geometry).

This derivative was introduced in \cite{cpr1} as the "spatial
derivative" of the orbit (in a more general context). There it was given
the expression
$$
\frac{D X}{dt}=\dot{X}+ \frac12 \left\{
[X,\kappa_{\gamma}(\dot{\gamma})]+[\dot{\gamma},\kappa_{\gamma}(X)]\right\}
$$
which coincides with our expression above after routine calculations.
\end{rem}

As already stated, $\um\subset {\cal P}(M_1)$ is a $C^\infty$ submanifold in the
norm
topology. In the strong topology (as metrized by $\|\ \|_2$) neither of these
spaces are manifolds. Nevertheless, with the trace Riemannian metric, ${\cal
P}(M_1)$
is well behaved: it is a complete metric space with both the Riemannian and the
$\|\ \|_2$-metric, which are equivalent, and any pair of projections can be
joined by minimal geodesics \cite{argrassmannians}. We do not know  if the
Riemannian metric and the $\|\ \|_2$-metric are equivalent in $\um$. One has
though the following inequality
$$
d_M(p_1,p_2)\ge \|p_1-p_2\|_2,
$$
where $p_1,p_2\in\um$ and $d_M$ denotes the Riemannian or geodesic distance.
This can be proved regarding $\um$ as a subset on $M_1$, and noting that any
curve in $\um$, as a curve in $M_1$, is longer than the straight line segment.

\medskip

Using Proposition \ref{completitud}, one can prove the following:

\begin{teo}\label{completo}
$\um$ is a complete metric space with the Riemannian metric $d_M$.
\end{teo}
\begin{proof}
Let $\{p_n\}$ be a Cauchy sequence in $\um$ for the metric $d_M$, with
$p_n=u_n\p u_n^*$. Then, by the above inequality, it is a Cauchy sequence for
the $2$-norm, and therefore there exists a projection $q$ such that
$\|p_n-q\|_2\to 0$. Therefore, by \ref{completitud}, $q=u\p u^*\in\um$.
Conjugating with $u$, we may suppose without loss of generality that $q=\p$.
Note that $\|p_n-\p\|_2\to 0$ is equivalent to $\|u_n-E(u_n)\|_2\to 0$, or
$\tau(E(u_n)^*E(u_n))\to 1$. Indeed
\begin{eqnarray}
\|p_n-\p\|_2^2 & =& \|u_n\p-\p u_n\|_2^2=\tau(\p)+\tau(u_n^*\p u_n)-\tau(\p u_n^*\p
u_n)-\tau(u_n^*\p u_n\p)\nonumber\\
&=&2\lambda-2\tau(\p u_n^*\p u_n\p)=2\lambda-2\tau(E(u_n^*)E(u_n)\p)=2\lambda(1-\tau(E(u_n^*)E(u_n))).\nonumber
\end{eqnarray}
Also note that, as in the proof of Proposition \ref{biyeccion}
$$
\|u_n-E(u_n)\|_2^2=1-\tau(E(u_n)^*E(u_n)).
$$
Let us apply now Lemma \ref{lemal2} to the elements $a_n=E(u_n)$ in the
algebra $N$, again, as in the proof of \ref{biyeccion}. Then there exist
unitaries $v_n$ in $N$ such that $\|E(u_n)-v_n\|_2\to 0$. It follows that
$\|u_n-v_n\|_2\to 0$, or equivalently,
$$
\|u_nv_n^*-1\|_2\to 0.
$$
In Proposition 4.4 of \cite{l2metric}, it was shown that in the unitary group of
a finite von Neumann algebra $M$ with trace $\tau$, the $2$-metric is equivalent
to the geodesic distance induced by the $2$-metric. Also in that paper, the
minimal geodesics were characterized, as the exponentials $\delta(t)=e^{itx}$
for $t\in[0,1]$, where $x^*=x$ and $\|x\|\le \pi$. It follows that there exist
$x_n\in M$ with $x_n^*=x_n$ and $\|x_n\|\le \pi$ which achieve the geodesic
distance between $1$ and $u_nv_n^*=e^{ix_n}$. And by the above facts, if
$\delta_n(t)=e^{itx_n}$, $t\in[0,1]$, then
$$
L_2(\delta_n)=\|x_n\|_2\to 0.
$$
Consider the curve
$$
\gamma_n(t)=e^{itx_n}\p e^{-itx_n} \in \um
$$
which join $\p$ to $u_nv_n^*\p v_n u_n^*=u_n\p u_n^*$. Then
$$
d_M(p_n,\p)\le L_2(\gamma_n)=\|x_n\p -\p x_n\|_2\le 2\|x_n\|_2\to 0.\qedhere
$$
\end{proof}

\begin{rem}
The geodesics of this connection that start at $q$ can be computed; they are the
curves of the form $\alpha(t)=e^{tz}q e^{-tz}$, with $z\in \h_q$.
\end{rem}
For any curve $\gamma\in \um$ (not necessarily a geodesic) with $\gamma(0)=q$,
there is an {\it horizontal lifting} $\Gamma$ to the unitary group $U_M$, which
is characterized by the following properties (see \cite{mr})
\begin{enumerate}
\item
$\ell_q(\Gamma)=\gamma$.
\item
$\Gamma(0)=1$.
\item
$\dot{\Gamma}\in \h_\gamma \Gamma $.
\end{enumerate}
Moreover, $\Gamma$ is also characterized as the unique solution of the linear
differential equation
\begin{equation}\label{levantadahorizontal}
\displaystyle\left\{ \begin{array}{l} \dot{\Gamma}=
\kappa_{\gamma}(\dot{\gamma})\Gamma \\ \Gamma(0)=1 \\  \end{array} \right.
\end{equation}
The following is an easy consequence of Lemma \ref{qiso} above.
\begin{prop}
Let $\gamma$ be a piecewise smooth curve in $\um$, and let $\Gamma$ be its
horizontal lifting. Then
$$
L_2(\gamma)=\sqrt{2\lambda} L_2(\Gamma).
$$
\end{prop}
\begin{proof}
Since $\Gamma q \Gamma^*=\gamma$, then
$$
L_2(\gamma)=\int_0^1 \|\dot{\gamma}\|_2 \ dt=
\int_0^1\|\dot{\Gamma}q\Gamma^*+\Gamma q \dot{\Gamma}^*\|_2 \ dt=\int_0^1
\|\Gamma(\Gamma^*\dot{\Gamma}q+q\dot{\Gamma}^*\Gamma)\Gamma^*\|_2 \ dt.
$$
Since $\Gamma$ is a curve of unitaries, differentiating $\Gamma^*\Gamma=1$ one
obtains that $\dot{\Gamma}^*\Gamma=-\Gamma^*\dot{\Gamma}$. Then the above
integral equals
$$
\int_0^1 \|\Gamma^*\dot{\Gamma} q-q \Gamma^*\dot{\Gamma}\|_2 \ dt=\int_0^1
\|\delta_q(\Gamma^*\dot{\Gamma})\|_2 \ dt=\sqrt{2\lambda} \int_0^1
\|\Gamma^*\dot{\Gamma}\|_2 \ dt=\sqrt{2\lambda}L_2(\Gamma).\qedhere
$$
\end{proof}

This observation implies that one can compare the lengths of the horizontal
liftings of the curves, instead of the curves themselves, and profit from the
geodesic structure of the unitary group. This is the point of view adopted in
this paper.

Let us finish this section proving that the horizontal lifting is the shortest
possible lifting for a curve in $\um$.
\begin{prop}
Let $\gamma(t)$, $t\in[0,1]$ be a piecewise smooth curve in $\um$. Suppose that
$u(t)\in U_{M}$ is a piecewise smooth lifting of $\gamma$. Then the horizontal
lifting $\Gamma$ of $\gamma$ is shorter than $u$:
$$
L_2(\Gamma)\le L_2(u).
$$
When using the usual operator norm, one has the estimate
$$
L_\infty(\Gamma)\le 2 L_\infty(u).
$$
\end{prop}
\begin{proof}
We may suppose without loss of generality that $\gamma$ starts at $\p$. Note
that $\gamma=\Gamma \p\Gamma^*=u\p u^*$. Then $u^*\Gamma\in \{\p \}'\cap M=N$.
Differentiating,
$$
\dot{u}^*\Gamma+u^*\dot{\Gamma}\in N.
$$
Note that $u^*\dot{\Gamma}=(u^*\Gamma) (\Gamma^*\dot{\Gamma})\in N\cdot\h\subset
\h$.
It follows that if $P_\h=I-E$ denotes the $\tau$-orthogonal projection onto
$\h$, one has
$$
P_\h(\dot{u}^*\Gamma)=-u^*\dot{\Gamma}.
$$
Therefore
$$
\|\dot{\Gamma}\|_2=\|-u^*\dot{\Gamma}\|_2=\|P_\h(\dot{u}^*\Gamma)\|_2\le
\|\dot{u}^*\Gamma\|_2=\|\dot{u}^*\|_2=\|\dot{u}\|_2.
$$
Thus
$$
L_2(\Gamma)=\int_0^1 \|\dot{\Gamma}\|_2 \ dt \le \|\dot{u}\|_2 \ dt=L_2(u).
$$
For the assertion corresponding $L_\infty$, we use the estimate
$$
\|P_\h(x)\|=\|x-E(x)\|\le 2 \|x\|.\qedhere
$$
\end{proof}

\begin{rem}
Using the argument above, one can prove that if $\gamma$ is a smooth curve in
$\um$, and $u$ is a {\it smooth} lifting of $\gamma$ such that
$L_2(u)=L_2(\Gamma)$, then there exists a unitary element $v_0$ commuting with
$\gamma(0)$ (for instance, if $\gamma(0)=\p$, then $v_0\in U_N$), such that
$$
u(t)=\Gamma(t) v_0 \ , \ \ t\in[0,1].
$$
That is, $u$ is essentially the horizontal lifting.
Indeed, in the proof above it is shown in fact that $\|\dot{\Gamma}(t)\|_2\le
\|\dot{u}(t)\|_2$. If $u$ is smooth, then $\|\dot{u}(t)\|_2$ is continuous in
the parameter $t$. Then the inequality above and the hypothesis that
$L_2(u)=\int_0^1\|\dot{u}(t)\|_2 d t=\int_0^1\|\dot{\Gamma}(t)\|_2 d
t=L_2(\Gamma)$ imply the equality
$$
\|\dot{\Gamma}(t)\|_2= \|\dot{u}(t)\|_2,
$$
for all $t$.
Therefore $\dot{u}^*\Gamma=P_\h(\dot{u}^*\Gamma)=-u^*\dot{\Gamma}$, or
equivalently
$$
0=\dot{u}^*\Gamma+u^*\dot{\Gamma}=\dot{u^*\Gamma},
$$
i.e. $u^*(t)\Gamma(t)=u^*(0)\Gamma(0)=v_0$ for all $t$, with $v_0$ commuting
with $\gamma(0)$.
\end{rem}

We may also compare the length of the horizontal lifting with the lengths of
liftings in $U_{M_1}$.
We shall need Lemma \ref{lemaseccion} from section 2.
\begin{prop}
Let $\omega(t)\in U_{M_1}$, $t\in[0,1]$, be a piecewise smooth lifting of
$\gamma(t)\in\um$.
Then
$$
L_2(\Gamma)\le \frac{1}{\sqrt{\lambda}} L_2(\omega).
$$
\end{prop}
\begin{proof}
By the lemma, there exists a piecewise smooth curve $u(t)\in U_M$ such that
$u(t)\p=\omega(t)\p$. In particular $u$ is a lifting (in $U_M$) of $\gamma$, and
therefore by the preceding proposition,
$$
L_2(\Gamma)\le L_2(u).
$$
On the other hand, since for $a\in M$ we have
$\|a\|_2=\frac{1}{\sqrt{\lambda}}\|a\p\|_2$, then
$$
L_2(u)=\frac{1}{\sqrt{\lambda}} L_2(u\p)=\frac{1}{\sqrt{\lambda}}
L_2(\omega\p)\le\frac{1}{\sqrt{\lambda}} L_2(u),
$$
were the last assertion follows from the fact that $\|x y\|_2\le \|x\|_2\|y\|$.
\end{proof}

\section{Convexity properties of the unitary group}

Our argument on the minimality of geodesics in Section 5 needs certain facts
concerning the
geometry of the unitary group $U_M$, which were proved in \cite{ar}. First
recall  that curves of the form $\delta(t)=ue^{itx}$ (with $x^*=x$ and $\|x\|\le
\pi$) have minimal length for the $k$-norms, $k\ge 2$ \cite{argrassmannians}.
Based on this fact, in Theorem 2.1 of \cite{ar} it was shown that if $F_2$
denotes the energy functional
$$
F_2(\gamma)=\int_0^1 \|\dot{\gamma}\|_2^2 dt=\int_0^1
\tau(\dot{\gamma}^*\dot{\gamma}) dt,
$$
for $\gamma$ a piecewise smooth curve in $U_M$, and $\gamma_s(t)$ is a  smooth
variation
of $\gamma$,
i.e.
$$
\gamma_s(t)\in U_M, s\in (-r,r) \ , \ \  t\in[0,1] \ , \ \  \gamma_0=\gamma,
$$
then the {\it first variation of the energy
functional} is
\begin{equation}\label{primeravariacion}
\frac12\frac{d}{d s}F_2(\gamma_s)|_{s=0}=\tau(x_0y_0)|_{t=0}^{t=1}-\int_0^1
\tau(\frac{d}{dt}[x_0] y_0) dt,
\end{equation}
where
$$
x_s(t)=\gamma_s(t)^* \frac{d}{dt}\gamma_s(t) \hbox{ and } y_s(t)=\gamma_s(t)^*
\frac{d}{ds}\gamma_s(t).
$$

The main result in \cite{ar} states the following:
\begin{prop}  \label{lemaclave}(Theorem 4.5 of \cite{ar})
Let $u_0$, $u_1$, $u_2\in U_M$, such that $\|u_i-u_j\|<\sqrt{2-\sqrt{2}}=r$.
Let $\delta(t)=u_1e^{tz}$ be the minimal geodesic joining $u_1$ and $u_2$.
Then $f(s)=d_k(u_0, \delta(s))^k$ ($d_k=$ geodesic distance induced by the
$k$-norm) is a convex function ($s\in[0,1]$), for $k$ an even integer.
\end{prop}

We shall use this result for the case $k=2$.

\bigskip

Take $z\in \h$; the curve $\alpha(t)=e^{tz}\p e^{-tz}$ is a geodesic of the
Levi-Civita connection of $\um$, joining $\p$ with $q=e^z\p
e^{-z}$. We have $\dot{\alpha}(t)=e^{tz}\left(z\p-\p
z\right)e^{-tz}$ hence
$$
\|\dot{\alpha}(t)\|_2=\|z\p-\p z\|_2=\|\dot{\alpha}(0)\|_2
$$
and
$$
L_2(\alpha)=\int_0^1\|\dot{\alpha}(t)\|_2dt=\sqrt{2\lambda}\|z\|_2
$$
by Lemma \ref{qiso} above.
The following estimates will prove useful. Recall that $\|\;\|$ denotes the
usual operator norm of the von Neumann algebra $M_1$:

\begin{lem}\label{uno}
For any $z\in\h=\ker E\cap M_{ah}$,
$$
\|\delta_{\p}(z)\| =\|z\p-\p z\| \ge
\sqrt{\lambda}\|z\|.
$$
\end{lem}
\begin{proof}
Since $\p z \p=0$,
$$
\|z\p-\p z\|^2 =\|(z\p-\p z)^2\| =\|z\p z+\p z^2\p\| =\|-z\p z-\p z^2\p\|
$$
Note that $-z\p z$ and $-\p z^2\p$ are both positive operators, and for $a,b\ge
0$, $\|a\| \le \|a+b\| $ hence
$$
\|z\p-\p z\|^2 \ge \|\p z^2\p\| =\|\p E(z^2)\| =\|E(z^2)\| \ge \lambda \|z\|^2
$$
by the definition of $E$ and the index properties.
\end{proof}

\begin{lem}
If $z\in\h$, with  $\|z\| <\sqrt{\lambda}$, then
$$
\|e^z\p e^{-z}
-\p\| \ge \|z\| (\sqrt{\lambda}-\|z\| ).
$$
\end{lem}
\begin{proof}
Note that $\sqrt{\lambda}<\pi$.
First, we rewrite the expression on the left:
\begin{eqnarray}
\|e^z\p e^{-z} - \p\|  &=&\|e^z\p-\p e^z\| =\|\delta_{\p}(z)+(e^z-1-z)\p - \p
(e^z-1-z)\|\nonumber\\
&\ge&\|\delta{_\p}(z)\| -\|(e^z-1-z)\p - \p (e^Z-1-z)\|\nonumber\\
& \ge& \|\delta{_\p}(z)\|
-2\|e^z-1-z\|\nonumber.
\end{eqnarray}
Now, since $\sigma(e^z-1-z)\subset \{ e^{it}-1-it:\mid t\mid \le \|z\|  \}$,
\begin{eqnarray}
\|e^z-1-z\| & =&\rho\left( e^z-1-z\right) =\sup\limits_{\mid t\mid \le\|z\|
}\sqrt{(\cos(t)-1)^2+(\sin(t)-t)^2}\nonumber \\
&=&\sup\limits_{\mid t\mid \le\|z\| }\sqrt{ 2-2\cos(t)-2t\sin(t)+t^2 }\nonumber\\
&=&\sqrt{2-2\cos\|z\| -2\|z\| \sin\|z\| +\|z\|^2  }.\nonumber
\end{eqnarray}
Using Taylor's series for the function under the square root, one can obtain the
bound $\frac12 \|z\| ^2$; hence
$$
\|e^z\p e^{-z} - \p\| \ge \|\delta_{\p}(z)\| -\|z\| ^2.
$$
This bound together with the one in the previous lemma gives the desired
inequality.
\end{proof}

Denote by $L_{\infty}(\gamma)$  the length of the curve $\gamma$,
measured in the  norm metric of  $\um$, namely
$$
L_{\infty}(\gamma)=\int_0^1\|\dot\gamma(t)\| dt.
$$

\begin{lem}\label{tres}
If $\gamma$ is a smooth curve in $\um$, and $\Gamma$ is the horizontal lift of
$\gamma$, then
$$
\|\Gamma(1)-1\|\le \frac{1}{\sqrt{\lambda}} L_{\infty}(\gamma).
$$
\end{lem}
\begin{proof}
First note that
$$
\|\Gamma(1)-1\| =\|\int_0^1 \dot{\Gamma(t)}dt\| \le \int_0^1 \|\dot{\Gamma(t)}\|
dt=\int_0^1\|\kappa_{\gamma}(\dot{\gamma})\| dt\le \int_0^1\|\kappa_{\gamma}\|
\, \|\dot\gamma\| dt.
$$
Now if $q\in\um$, $\|\kappa_q\| =\|\kappa_{\p}\| $ because the action of the
unitary group is isometric. Since $\kappa_{\p}=\delta_{\p}^{-1}$, where
$\delta_{\p}:\h\to T_{\p}\um$, by the bound in lemma \ref{uno}, we have
$\|\kappa_{\p}(v)\| \le \frac{1}{\sqrt{\lambda}}\|v\| $, namely $\|\kappa_q\|
\le\frac{1}{\sqrt{\lambda}}$ for any $q\in \um$. Hence
$$
\|\Gamma(1)-1\| \le \frac{1}{\sqrt{\lambda}}\int_0^1\|\dot{\gamma(t)}\|  dt=
\frac{1}{\sqrt{\lambda}}L_{\infty}(\gamma).\qedhere
$$
\end{proof}

We say that a subset $V\subset U_M$ is {\it geodesically convex} if for any pair
of elements $u_0,u_1\in V$ with $\|u_0-u_1\|<2$, the unique minimal geodesic
of $U_M$ joining them lies in $V$.

\begin{teo}\label{corto}
There exists a positive constant $R>0$ with the following property. Suppose
that $q_0,q_1\in\um$ with
 $q_1=e^z q_0 e^{-z}\in \um$ for $z\in \h_{q_0}$ with $\|z\| \le R$, and let
$\alpha(t)=e^{tz} q_0 e^{-tz}$ be the geodesic joining $q_0$ and $q_1$. If
$\gamma$ is any other smooth curve joining $q_0$ and
$q_1$, then either
$$
L_{\infty}(\gamma)\ge L_{\infty}(\alpha),
$$
or
$$
L_2(\gamma)\ge L_2(\alpha).
$$
\end{teo}
\begin{proof}
Clearly all the norms involved are unitarily invariant, therefore we may suppose
$q_0=\p$ without loss of generality.
Let $\Gamma$ be the unitary lift of $\gamma$ in $U_M$. If $L_\infty(\gamma)\le
L_\infty(\alpha)$ , by Lemma \ref{tres} above
$$
\|\Gamma(1)-1\|  \le \frac{1}{\sqrt{\lambda}}\|z\p-\p z\|\le
\frac{2}{\sqrt{\lambda}}\|z\|.
$$
Therefore we may adjust $R$ in order that  $1,e^z$ and $\Gamma(1)$ stand in the
situation where Proposition \ref{lemaclave} applies, namely, that they lie
closer than $\sqrt{2-\sqrt{2}}$ in norm.
Let $\delta(s)=e^ze^{sw}$ be the unique minimizing geodesic joining
$\delta(0)=e^z$ and $\delta(1)=e^ze^w=\Gamma(1)$ in $U_M$ (which is minimizing
for all $k$-norms, and therefore in particular for the functional $L_2$
\cite{argrassmannians}). Then by Proposition \ref{lemaclave}, the map
$$
f(s)=d_2^2(1, \delta(s)), \ \ s\in[0,1]
$$
is convex. We claim that $f'(0)=0$, so that $f$ has an absolute minimum at
$s=0$.
Indeed, note that
$$
\|\delta(s)-1\|=\|1-e^{sw}\|\le \|1-e^w\|=\|1-\Gamma(1)\|<2.
$$
Therefore the antihermitic logarithm
$$
log:\{u\in U_M:\|u-1\|<2\}\to \{x\in M_{ah}: \|x\|<\pi\}
$$
is well defined.
Let $\gamma_s(t)=e^{t log(\delta(s))}$.
Note that $\gamma_s(t)$ is a smooth variation of $\gamma_0=\alpha$. Also note
that at
each $s$ it consists of minimizing geodesics, because
$\|log(\delta(s))\|<\pi$. Then
$$
f(s)=L_2(\gamma_s)^2=\|log(\delta(s))\|_2^2.
$$
Note also that $f(s)=F_2(\gamma_s)$. Then $f'(0)$ can be
computed using the first variation formula (\ref{primeravariacion}). In our case
$$
x_s=\gamma_s^*\frac{d}{dt}\gamma_s(t)=log(\delta(s))
$$
is independent of $t$, and therefore (\ref{primeravariacion}) reduces to
$$
f'(0)=2\tau(z y_0(1))-\tau(zy_0(0)).
$$
Note that $\Gamma(1)$ and $e^z$ lie in the fibre of $q_1$ (both
$\gamma=\Gamma\p\Gamma^*$ and $\alpha$ have the same endpoints), which is of the
form $e^z U_N$. Clearly this set is geodesically convex \cite{ar}, implying that
$\delta(s)$ lies in this fibre, and therefore $\dot{\delta}(0)\in e^z N_{ah}$.
At
$t=0$, $\gamma_s(0)=1$ for all $s$, therefore $y_s(0)=0$. At $t=1$,
$\gamma_s(1)=\delta(s)$,
so that $\gamma_s(1)=\delta^*(s)\dot{\delta}(s)$. Then
$y_0(1)=\delta^*(0)\dot{\delta}(0)=e^{-z}\dot{\delta}(0)\in N_{ah}$. Then
$\tau(z y_0(0))=0$ and
$$
\tau(z y_0(1))=\tau(E(z y_0(1)))=\tau(E(z) y_0(1))=0
$$
because $E(z)=0$.

\medskip

Our claim proved, it implies that if we denote $A(t)=e^{tz}$, and $B(t)=e^{t
\log(\Gamma(1))}$, $t\in[0,1]$, then
$$
L_2(A)\le L_2(B).
$$
On the other hand, by the fact on minimality of curves in
$U_M$ \cite{argrassmannians}, one has that $L_2(B)\le L_2(\Gamma)$. Therefore
$$
L_2(A)\le L_2(\Gamma).
$$
Multiplying both members by $\sqrt{2\lambda}$ one gets
$$
L_2(\alpha)=\sqrt{2\lambda}L_2(A)\le
\sqrt{2\lambda}L_2(\Gamma)=L_2(\gamma).\qedhere
$$
\end{proof}

\begin{rem}\label{Rsub1}
Note that there exists a constant $R_1>0$ such that $\|q_1-q_0\|<R_1$ is
equivalent to the existence of $z\in \h_{q_0}$ with $e^zq_0e^{-z}=q_1$ and
$\|z\|<R$. Indeed, for the norm topology, due to the differential structure
of $\um$, the map
$$
\exp_{q_0}:\h_{q_0}\to \um,\; \exp_{q_0}(z)=e^zq_0e^{-z}
$$
is a local diffeomorphism by the inverse function theorem.

Therefore the result above can be rephrased replacing the requirement
$\|z\|<R$ by $\|q_0-q_1\|<R_1$.
\end{rem}

The theorem above states that geodesics are short among curves which are {\it a
priori} short in the norm Finsler metric. In the next section we show a more
general minimality result, but let us state this sentence with more precision,
for we shall need it later.
\begin{coro}\label{Rsub2}
There exists a numbre $R_2>0$ such that if $\gamma$ is a smooth curve with
$L_\infty(\gamma)\le R_2$, then there exists a geodesic $\delta$ joining the
same endpoints as $\gamma$, with $L_2(\delta)\le L_2(\gamma)$.
\end{coro}
\begin{proof}
Choose $R_2<R_1$ in order to assure the existence of $\delta$: if $q,r$ are
the endpoints of $\gamma$, then $\|q-r\|\le L_\infty(\gamma)\le R_2<R_1$ and
Remark \ref{Rsub1} applies. Then the proof  follows as in the theorem above,
adjusting  $R_2$ additionally in order that
$$
\|\Gamma(1)-1\|\le \frac{1}{\sqrt{\lambda}}L_\infty(\gamma)\le
\frac{R_2}{\sqrt{\lambda}}<\sqrt{2-\sqrt{2}}.\qedhere
$$
\end{proof}

\section{Geodesics as unique minimal curves}

A continuous curve $\alpha$ is a \textit{piecewise smooth geodesic} if it
consists of a finite collection $\{\alpha_i\}$ of geodesic arcs glued together,
or in other words, a polygonal path with geodesic edges. Let us start with a
brief result:

\begin{prop}
Let $\gamma\subset\um$ be any smooth curve joining $q$ to $r$. Then there exists
a continuous
 piecewise smooth geodesic $\alpha=\cup \alpha_i$ which joins $q$ and $r$ such
that
$$
\sum L_2(\alpha_i)\le L_2(\gamma).
$$
\end{prop}
\begin{proof}
Suppose $\gamma$ is parametrized in the interval $[0,1]$. Then there exists a
partition $0=t_0<t_1<\dots<t_n=1$ of $[0,1]$ such that
$\|\gamma(t_i)-\gamma(t_{i-1})\|<R_1$ for $i=1\dots n$ with $R_1$ as in
Remark \ref{Rsub1}. Then there exist geodesics $\alpha_i$, minimizing for the
functional length $L_2$, which join $\gamma(t_{i-1})$ with $\gamma(t_i)$. Then
clearly $\alpha=\cup \alpha_i$ is shorter for the 2-metric than $\gamma$.
\end{proof}

Our main result on uniqueness of geodesics as minimal curves follows.
\begin{teo}\label{unici}
Let $\gamma\subset\um$ be a piecewise smooth curve which is short (i.e
minimizing) for the 2-metric. Then $\gamma$ is a  geodesic of the Riemannian
connection of $\um$.
\end{teo}
\begin{proof}
By the previous lemma, there exists a piecewise smooth geodesic $\alpha=\cup
\alpha_i$ which is shorter than $\gamma$. Then  $\alpha$ is also minimizing (in
fact, $L_2(\gamma)=\sum L_2(\alpha_i)=L_2(\alpha))$. Moreover, the elements
$t_0,t_1,...,t_n$ of the partition provide points $\gamma(t_i)$ which lie both
in $\gamma$ and $\alpha$. We claim that $\alpha$ is smooth, i.e. a geodesic.
This proves our result, since the partition can be arbitrarily refined to
contain as many (finite) points in common between $\gamma$ and $\alpha$, and the
smoothness of $\alpha$ proves that all these polygonals are in fact the same
geodesic.

Asumme that $\alpha$ is not smooth to arrive to a contradiction. Namely assume
there is at least one point $q_i$ where $\dot\alpha_i(1^-)\ne
\dot\alpha_{i+1}(0^+)$. Since $\alpha_i(t)=e^{tz^+}q_ie^{-tz^+}$ and
$\alpha_{i+1}(t)=e^{tz^-}q_ie^{-tz^-}$, then
$$
\Delta\dot\alpha_i=\dot\alpha_{i+1}(0^+)- \dot\alpha_i(1^-)=[z^+-z^-,q_i]
$$
is a nonzero vector in $\h_q$. We may choose a variation $\gamma_s$ of
$\alpha=\gamma_0$ which is constantly identical to $\alpha$ except in a
neighbourhood of $q_i$, and such that the variation field
$V(t)=\frac{\partial}{\partial s}\mid_{s=0}\gamma_s$ equals this vector in
$t=1$, namely $V(1)=\Delta\dot\alpha_i$.

According to the classic first variation formula in a Riemannian manifold (cf.
\cite{lee} for example),
$$
\frac{d}{ds}|_{s=0}L_2(\gamma_s)=-\int_0^1 \pei V, D_t\dot\gamma \ped dt - \sum
\pei V(t_i),\Delta_i\dot\gamma \ped.
$$
Since $\alpha$ consists of piecewise geodesics and $\alpha$ is also a critical
point of the length distance (it is minimizing for the 2-metric) we obtain
$$
0=\frac{d}{ds}\mid_{s=0}L_2(\gamma_s)=\|[z^+-z^-,q_i]\|_2
$$
which is absurd.
\end{proof}

Our main result concerning existence of minimal geodesics follows.
\begin{teo}
Let $\gamma$ be a smooth curve starting at $\p$ that stays in a neighbourhood of
$\p$ of radius $R_2/2$. Then $\gamma$ is longer than the geodesic $\delta$ that
joins its endpoints, namely $L_2(\gamma)\ge L_2(\delta)$.
\end{teo}
\begin{proof}
We partition $\gamma$ in the same fashion as before, making sure that each piece
is shorter than $R_2$ in the $\infty$-metric. Now we replace each piece of
$\gamma$ with a minimizing geodesic; let us call the breaking points $q_i$,
where $q_0=\p$ and $q_n$ is the endpoint of $\gamma$. Clearly the first geodesic
is shorter than the first piece of $\gamma$; note that by the triangle
inequality, distance between $q_2$ and $q_0=\p$ is shorter than $R_2$, hence we
may replace the segments $q_0q_1$, $q_1q_2$ with a shorter geodesic
$\delta_{11}$. This geodesic is shorter than the segment of $\gamma$ which joins
$q_0$ to $q_2$, and if we proceed inductively, we end up with a smooth geodesic
$\delta$ which joins $q_0=\p$ to $q_n$, and this geodesic $\delta$ is clearly
shorter than $\gamma$.
\end{proof}

\section{The tangent and normal bundles} 
Let us finish this paper  characterizing the tangent spaces of $\um$ as subspaces of the tangent bundle of $\grass$, and thereon give a characterization of the tangent bundle of $\grass$ and the normal space of $\um$ at $\p$, using only elements in the algebra $M$. We show in this description that curvature of $\um$ in $\grass$ is  related to the  set $\{ x^*=-x \in \ker E:\; x^2\in N\}\subset M$.

\begin{prop}
Let $q\in\um$, then
$$
(T\um)_q=\{y\in (T\grass)_q: E_1(y)=0\}.
$$
\end{prop}
\begin{proof}
Note  that $\um$ consists of projections $q$ in $M_1$ such that $E_1(q)=\lambda$
(see \cite{popa}). Therefore, if $q(t)$ is a curve in $\um$, then
$E_1(\dot{q}(t))=0$. In ohter words, $(T\um)_q\subset \{y\in (T\grass)_q:
E_1(y)=0\}$. Conversely, suppose that $y\in (T\grass)_q$ and $E_1(y)=0$. If
$q=u\p u^*$, by conjugating with $u^*$ one may suppose $q=\p$. Then there exists
$x\in M_1$, $x^*=-x$ such that $y=x\p-\p x$. By the  properties of the basic
construction, there exists $m\in M$ (a priori not necessarily anti-hermitics)
such that $x\p=m\p$. Then $-\p x=(\p x)^*=\p m^*$. Then
$$
0=E_1(y)=E_1(m\p+\p m^*)=\lambda m+\lambda m^*,
$$
i.e. $m\in M_{ah}$, and $y=m\p - \p m \in (T\um)_\p$.
\end{proof}

Now recall that $(T{\cal O}(\p))_{\p}=\{z\p-\p z : z^*=-z, z\in \ker E\}$, where $\ker E=N^{\perp}\subset M$.
\begin{prop}We have
$$
(T\grass)_{\p}=\{x\p+\p x^*:x\in N^{\perp}\}\;\mbox{ and }\;
(T\um)_{\p}^{\perp}=\{x\p+\p x : x^*=x,\; x\in N^{\perp}\}.
$$
In particular $\displaystyle (T \grass)_{\p}\simeq N^{\perp}\simeq N_{ah}^{\perp}\oplus N_{h}^{\perp}\simeq (T\um)_{\p}\oplus (T\um)_{\p}^{\perp}.$

\medskip

The Riemannian exponential map of $\grass$ is given by 
$$
\exp_{\p}^{\grass}(x\p+\p x^*)=e^{(x\p-\p x^*)}\p e^{-(x\p-\p x^*)}.
$$

In the decomposition relative to Jones' projection $\p$ we have
$$
e^{x\p -\p x^*}=\left( \begin{array}{ccc} \cos\left(\sqrt{E(\mid x\mid^2})\right)\p  & & - \mbox{\rm sinc}\left(\sqrt{E(\mid x\mid^2})\right) \p x^* \\ \\ x\p \mbox{\rm sinc}\left(\sqrt{E(\mid x\mid^2})\right) & & \cos(\sqrt{x\p x^*}) \end{array}   \right).
$$
where $\mbox{\rm sinc}(z)=\sin(z)z^{-1}$. Moreover, 
$$
e^{(x\p -\p x^*)}\p e^{-(x\p -\p x^*)}=\left( \begin{array}{ccc} \cos^2\left(\sqrt{E(\mid x\mid^2})\right)\p  & &  \mbox{\rm sinc}\left(2\sqrt{E(\mid x\mid^2})\right) \p x^* \\ \\ x\p \mbox{\rm sinc}\left(2\sqrt{E(\mid x\mid^2})\right) & & \sin^2(\sqrt{x\p x^*}) \end{array}   \right).
$$
The action $\p\mapsto e^{(x\p -\p x^*)}\p e^{-(x\p -\p x^*)}$ is effective if $\|x\|<\pi$.
\end{prop}
\begin{proof}
Let $v\in (T\grass)_{\p}=\{w\p -\p w:w\in M_1, w \mbox{ is } \p -\mbox{codiagonal}\}$. Then  we have $v=w\p-\p w=R(w)\p+\p R(w)^*$ with $R(w)\in M$. Now $E(R(w)) \p=\p R(w)\p=\p w\p =0$, hence $R(w)\in\ker E=N^{\perp}$. On the other hand, for given $x\p+\p x^*$ with $x\in \ker E$, let $w=x\p -\p x^*$.  Then $w=w\p+\p w$ ($w$ is $\p$-codiagonal) and $w\p-\p w=x\p +\p x^*$. The isomorphism and the formula for the Riemannian exponential are now clear from the results on the previous sections.

The formula for the exponential $e^{x\p -\p x}$ can be deduced from the elementary formula
$$
\exp\left(  \begin{array}{ccc} 0 & & -Y^* \\ \\ Y & & 0 \end{array}   \right)= \left(  \begin{array}{ccc} \cos\mid Y\mid & & -\mbox{\rm sinc}\mid Y^*\mid Y \\ \\ Y\mbox{\rm sinc}\mid Y^*\mid & & \cos\mid Y^*\mid \end{array}   \right)
$$
putting $Y=x\p$ and recalling that $\p x\p=E(x)\p =\p E(x)$.

Note that, if $e^{x\p-\p x^*~}$ commutes with $\p$, then $\cos^2\left(\sqrt{E(\mid x\mid^2})\right)\p=\p$ (equivalently $\cos^2\left(\sqrt{E(\mid x\mid^2)}\right)=1$, since $E(\mid x\mid^2)\in N$). This equation proves that $\sigma(   \sqrt{ E( \mid x\mid^2)})\subset \{k\pi\}_{k\in \mathbb Z}$. Since $\|E(\mid x\mid^2)\|\le \|x\|^2<\pi^2$, we obtain  $\sigma(\sqrt{ E( \mid x\mid^2)})=\{0\}$. Since $\sqrt{ E( \mid x\mid^2)}$ is a positive operator, it must be zero. Equivalently, since $E( \mid x\mid^2 )\ge \lambda \mid x\mid^2$, we obtain $x=0$.
\end{proof}

\smallskip

\begin{teo}
The geodesic of $\grass$ given by $\gamma_x(t)=e^{t(x\p -\p x^*)}\p e^{-t(x\p -\p x^*)}$ is a geodesic of ${\cal O}(\p)$ iff $x^*=-x$ and $x^2\in N$, and in this case we have 
$$
\gamma_x(t)=e^{tx}\p e^{-tx}=\p \cos^2(t\mid x\mid)+ u\p u^*\sin^2(t\mid x\mid) +\frac12 [u,\p ]\sin(2t\mid x\mid).
$$ 
where $x=u\mid x\mid$, with $u^*=-u\in N^{\perp}$ a partial isometry such that $u^*u=uu^*=-u^2\in N$, and $\mid x\mid=\sqrt{-x^2}\in N$.
\end{teo}
\begin{proof}
Assume first that $x^*=-x$ and $x^2\in N$. Let $x=u\mid x\mid$ be the polar decomposition of $x$, note that $u^*=-u$ commutes with $x$. Then $e^{tx}\p=e^{-i(ix)}\p=\cos(ix)-i\sin(ix)=\cos(\mid x \mid)+ u\sin (\mid x\mid)$, because $(ix)^{2k}=i^{2k} x^{2k}=(-1)^k u^{2k}\mid x\mid^{2k}=\mid x\mid^{2k}$. Likewise, $(ix)^{2k+1}=u i \mid x\mid^{2k+1}$. On the other hand, from the previous proposition we obtain
$$
e^{x\p -\p x*}\p=e^{x\p+\p x}\p=\left(\cos\mid x\mid+x\, {\rm sinc}\mid x\mid\right)\p=\left(\cos\mid x\mid+u\sin\mid x\mid\right)\p,
$$
because $\sqrt{E(\mid x\mid^2)}=\mid x\mid$ since $\mid x\mid^2\in N$. Hence 
$e^{t(x\p -\p x^*)}\p =e^{tx}\p$ in this case.

Assume now that $\gamma_x(t)=e^{t(x\p -\p x^*)}\p e^{-t(x\p -\p x^*)}\p=e^{z}\p e^{-tz}\in {\cal O}(\p)$ with $x\in N^{\perp}$. Computing the first derivative at $t=0$ and multiplying by $\p$ on the right we obtain $x=z$. Computing the second derivative at $t=0$ we obtain $E(x^2)-x^2=0$, hence $x^2\in N$.

If $x=u\mid x\mid$, we have $u^*x=u^*u\mid x\mid=\mid x\mid$. Note that $x^2\in N$ iff $\mid x\mid\in N$. We can write $u=x h(\mid x\mid)$ with $h$ a Borel function (\cite{simon}  ex. IV.3). This shows that $u\in N^{\perp}$. Since $x^*=-x$, we have $u^*=-u$, and $u^*u=u^*x h(\mid x\mid)=\mid x\mid h(\mid x\mid)\in N$. The formula for the geodesic now follows from the proposition above.
\end{proof}

\begin{coro}
If $\lambda=\frac12$ (i.e. $[N:M]=2$), then ${\cal O}(\p)$ is a totally geodesic submanifold of $\grass$.
\end{coro}
\begin{proof}
If  $\lambda=\frac12$, the extension $[N:M]$ can be represented with $2\times 2$ matrices with entries in $N$, namely
$$
M=\left( \begin{array}{cc}n & \theta(n') \\ n'&  n \end{array} \right) \quad \mbox{and}\quad N=\left( \begin{array}{ccc}n & 0  \\ 0 & n \end{array} \right)
$$
where $\theta$ is an isomorphism of $N$ of order 2. Then $N^{\perp}=\left( \begin{array}{ccc} 0 & \theta(n')  \\ n' & 0 \end{array} \right)$, hence $(N^{\perp})^2\subset N$.

\end{proof}

\vskip1cm

{
\noindent
Esteban Andruchow and Gabriel Larotonda\\
Instituto de Ciencias \\
Universidad Nacional de General Sarmiento \\
J. M. Gutierrez 1150 \\
(1613) Los Polvorines \\
Argentina  \\
e-mails: eandruch@ungs.edu.ar, glaroton@ungs.edu.ar
}

\end{document}